\documentclass[11pt]{amsart}

\usepackage{amsmath,amssymb}
\usepackage{graphicx}
\usepackage{import}
\usepackage{amscd}
\usepackage{wrapfig}
\usepackage{fullpage}
\usepackage{epsfig}
\usepackage{mathptmx}
\numberwithin{equation}{section}

\newtheorem{theorem}{Theorem}[section]

\theoremstyle{definition}

\theoremstyle{remark}

\newcommand{\C}{{\mathbb C}}

\newcommand{\R}{{\mathbb R}}

\newcommand{\T}{{\mathcal T}}
\renewcommand{\H}{{\mathbb H}}

\newcommand{\x}{{\mathbf{x}}}
\newcommand{\y}{{\mathbf{y}}}

\newcommand{\J}{{\mathcal J}}

\newcommand{\mcg}{\mathcal{MC\kern0.04emG}}
\newcommand{\ml}{\mathcal{M \kern0.07emL}}
\newcommand{\pmf}{\mathcal{P \kern0.07em M \kern0.07em F}}

\newcommand{\teichmuller}{Teichm{\"u}ller{ }}

\makeatletter
 \let\c@theorem=\c@subsection
 \let\c@conjecture=\c@subsection
 \let\c@lemma=\c@subsection
 \let\c@proposition=\c@subsection
 \let\c@claim=\c@subsection
 \let\c@question=\c@subsection
 \let\c@criterion=\c@subsection
 \let\c@vfconj=\c@subsection
 \let\c@definition=\c@subsection
 \let\c@notation=\c@subsection
 \let\c@remark=\c@subsection
 \let\c@example=\c@subsection
 \let\c@equation=\c@subsection
 \let\c@figure=\c@subsection
 \let\c@wrapfigure=\c@subsection

\makeatother

\begin{document}

\title[linear inv]{The limit set of the handlebody set has measure zero.}
\author[Gadre]{Vaibhav Gadre}
\address{Mathematics, Harvard University, One Oxford Street, Cambridge,  MA 02138, USA}
\email{vaibhav@math.harvard.edu}

\begin{abstract}
This note fixes a small gap in the proof in \cite{Ker2} that the limit set of the handlebody set has measure zero.
\end{abstract}

\maketitle

%%%%%%%%%%%%%%%%%%%%%%%%%%%%%%%%%%%%
\section{Introduction}
%%%%%%%%%%%%%%%%%%%%%%%%%%%%%%%%%%%%
For an orientable surface $S$ of genus $g$, the mapping class group $\mcg(S)$ is the group of orientation preserving diffeomorphisms of $S$ modulo isotopy. The \teichmuller space $T(S)$ is the space of marked complex structures on $S$, and is homeomorphic to an open ball in $\R^{6g-6}$. Thurston showed that $T(S)$ is naturally compactified by $\pmf(S)$, the space of projective measured foilations on $S$. The space $\pmf(S)$ is homeomorphic to $S^{6g-7}$. 

In many respects, the action of $\mcg(S)$ on the compactification $\overline{T(S)}$ is similar to the action of a Kleinian group on $\overline{\H^3} = \H^3 \sqcup \C \mathbb{P}^1$. There is a natural Lebesgue measure class $\ell$ on $\pmf(S)$ that is ergodic for this action. In a chart given by a complete train track, the natural volume form on the space of normalized weights carried by the track is a measure in this class. For finitely generated Kleinian groups, the Ahlfors measure conjecture, states that the limit set in $\C \mathbb{P}^1$ has either zero or full measure. See \cite{Cal}. It is natural to wonder if a similar property is true for limit sets in $\pmf(S)$ of subgroups of $\mcg(S)$. A handlebody $H$ gives an interesting subgroup of $\mcg(S)$ called the handlebody group $G_H$. In \cite{Ker2}, Kerckhoff gave an elegant proof that the limit set $\Lambda(G_H)$ of $G_H$ has measure zero. However, it uses an earlier result from \cite{Ker1} which states that almost every splitting sequence of complete train tracks gets uniformly distorted infinitely often. Subsequently, the proof of this earlier result was discovered to be incomplete. 

We do not know how to prove uniform distortion for all complete train tracks. However, we show in \cite{Gad} that uniform distortion does hold for certain expansions of complete train tracks with a {\em single switch}. As analogs of interval exchange maps, we call such tracks {\em complete non-classical exchanges}. Just as for interval exchange maps, Rauzy induction is well defined for such tracks, and expansions are given by iterated Rauzy induction. See \cite{Gad} for the details. 

In this note, we show that the uniform distortion result from \cite{Gad} is sufficient to complete the proof in \cite{Ker2} that $\ell(\Lambda(G_H)) = 0$.

\subsection{The limit set of the handlebody set:}
Let $H$ be a handlebody with boundary surface $S$, an orientable surface of genus $g$. The handlebody group $G_H$ is the subgroup of $\mcg(S)$ consisting of classes that have a representative that extends over $H$. The limit set $\Lambda(G_H)$ is the smallest non-empty closed invariant subset for the action of $G_H$ on $\pmf(S)$. There is a natural Lebesgue measure class $\ell$ on $\pmf(S)$. The main theorem here is:

\begin{theorem}\label{limit-set}
For any handlebody $H$, 
\[
\ell\left( \Lambda(G_H) \right) = 0.
\]
\end{theorem}
 
%%%%%%%%%%%%%%%%%%%%%%%%%%%%%%%%%%%%%%%%%%%%%
\section{Masur's description of the limit set}
%%%%%%%%%%%%%%%%%%%%%%%%%%%%%%%%%%%%%%%%%%%%% 
Following the notation in \cite{Ker2}, let $\mathcal{B}$ be the set of isotopy classes of essential simple closed curves on $S$ that bound discs in $H$. A {\em cut system} is a collection of simple closed curves $C= \{C_1, \cdots , C_g\}, C_i \in \mathcal{B}$ which together with the disks they bound, cut $H$ into a 3-ball. As Masur proved in \cite{Mas}, a curve $\gamma \in \mathcal{B}$ if and only if for every cut system $C$, after $\gamma$ is isotoped to have minimal intersections with $C$, one of the following conditions is satisfied:
\begin{enumerate}
\item For every $C_i$, the intersection $C_i \cap \gamma = \varnothing$.
\item For some $C_i$, the curve $\gamma$ has {\em returning arcs} to $C_i$ i.e., $\gamma$ intersects $C_i$ and then before intersecting any $C_j$ (which includes $j=i$), it returns from the same side it just left to intersect $C_i$ again.
\end{enumerate}
Kerckhoff shows that after passing to $\Lambda(G_H)$, the above conditions persist i.e., the foliations in $\Lambda(G_H)$ satisfy either (1) or (2). The set of foliations which satisfy (1) is directly seen to be a measure zero set. So it suffices to show that the subset $\mathcal{R}_H$ of $\Lambda(G_H)$ consisting of foliations that have returning arcs for all cut systems, has measure zero.

%%%%%%%%%%%%%%%%%%%%%%%%%%%%%%%%%%%%%%%%%%%%%%
\section{Train tracks}
%%%%%%%%%%%%%%%%%%%%%%%%%%%%%%%%%%%%%%%%%%%%%%
Roughly speaking, a train track on $S$ is {\em complete} if all its complementary regions in $S$ are ideal triangles. Technically, the definition requires that the track be recurrent and transverse recurrent, but this will always be the case for the tracks considered here, so we skip the details and refer the reader to \cite{Pen}. For a train track $\tau$, let $P(\tau)$ denote the set of projective measured foliations carried by $\tau$. An assignment of non-negative weights to the branches of $\tau$ is said to be carried by $\tau$ if at all switches of $\tau$, the sum of weights for incoming branches is equal to the sum of the weights for outgoing branches. The set $P(\tau)$ can be identified with the set of normalized weights carried by $\tau$, where the normalization is that the sum of the weights is 1. When $\tau$ is complete, the set $P(\tau)$ gives a chart on $\pmf(S)$. The volume form on the set of normalized weights carried by $\tau$ defines a measure $\ell_\tau$ on $P(\tau)$. The transition maps for these charts are piecewise linear invertible maps, and the derivatives of the transition maps are not bounded. This means that only the Lebesgue measure class $\ell$ is defined by this process, and care is necessary while considering an infinite sequence of train track charts.

In Theorem 1 of \cite{Ker2}, Kerckhoff proves the following key result for train track charts:

\begin{theorem}\label{proportion}
There is a constant $0< K< 1$ such that for any complete train track $\tau$, there is subset $P(\tau, H) \subset P(\tau)$ such that 
\[
\frac{\ell_\tau(P(\tau, H))}{\ell_\tau(P(\tau))} > K
\]
and $P(\tau, H) \cap \mathcal{R}_H = \varnothing$. 
\end{theorem}
\noindent In other words, for any complete train track $\tau$, a definite proportion of $P(\tau)$ is disjoint from $\mathcal{R}_H$, where the proportion is calculated in the measure $\ell_\tau$. As a particular case, the above theorem holds for complete non-classical exchanges.
 
The proof in \cite{Ker2} then proceeds as follows: Cover $\pmf(S)$ by a finite collection of train track charts. It suffices to prove that $\mathcal{R}_H$ has measure zero in each of these charts. For any complete train track $\tau$ in this collection, Theorem~\ref{proportion} states that a definite proportion of it is disjoint from $\mathcal{R}_H$. Then, $\tau$ is split enough number of times such that the complement is a union of stages in splitting expansions. Applying Theorem~\ref{proportion} to each of these stages yields a definite proportion of each of them disjoint from $\mathcal{R}_H$. However, this proportion is being measured separately in each of the stages and not using the measure $\ell_\tau$. 

To take care of this issue, Kerckhoff uses the uniform distortion result from \cite{Ker1} (whose proof is incomplete). The uniform distortion result implies that up to leaving out a set of $\ell_\tau$-measure zero, the complement can be split into a (possibly infinite) union of uniformly distorted stages. If a stage is uniformly distorted, the proportion that is calculated in its measure, changes by a universally bounded amount when calculated using $\ell_\tau$. Continuing to split what remains into uniformly distorted stages and iterating Theorem~\ref{proportion} proves that $\mathcal{R}_H$ has measure zero in the initial chart $P(\tau)$. 

In the next section, we will state and explain the uniform distortion theorem from \cite{Gad} and use that instead in this final part of Kerckhoff's argument.

%%%%%%%%%%%%%%%%%%%%%%%%%%%%%%%%%%%%%%%%%%%%%
\section{Uniform distortion and proof of Theorem~\ref{limit-set}}
%%%%%%%%%%%%%%%%%%%%%%%%%%%%%%%%%%%%%%%%%%%%%

\subsection{Complete non-classical exchanges:} 
Complete train tracks $\tau$ with a {\em single switch} are called {\em complete non-classical exchanges}. Because of the single switch condition, these are analogs of interval exchange maps. Just as with interval exchange maps, Rauzy induction can be defined for non-classical exchanges. Iterated Rauzy induction associates an expansion to a measured foliation in $P(\tau)$. As expected, almost every measured foliation in $P(\tau)$ has an infinite expansion. 

Similar to interval exchange maps, a finite Rauzy sequence giving $\tau'$ from $\tau$ is encoded by an integer matrix $Q$ with non-negative entries and determinant 1. This gives a projective linear map $\T Q$ from the space $P(\tau')$ of normalized weights carried by $\tau'$ to the space $P(\tau)$ of normalized widths associated to $\tau$. A finite Rauzy sequence is $C$-uniformly distorted if the Jacobian of its projective linear map $\J(\T Q)$ satisfies
\[
\frac{1}{C} < \frac{\J(\T Q) (\x) }{\J(\T Q)(\y)} < C
\]
for all points $\x, \y \in P(\tau')$. In \cite{Gad}, we prove the following theorem:
\begin{theorem} [Uniform Distortion]\label{uniform-distortion}
There is a universal constant $C>1$ that depends only on the genus $g$, such that for almost every projective measured foliation in $P(\tau)$, its Rauzy expansion gets $C$-uniformly distorted infinitely often. 
\end{theorem}
$C$-uniformly distorted stages are important because up to a universal constant, relative probabilities remain unchanged. To be precise, there is a universal constant $0< c< 1$ depending only on the genus $g$, such that for any $C$-uniformly distorted stage with exchange $\sigma$ and associated projective linear map $\T Q: P(\sigma) \to P(\tau)$, and for any measurable subset $A \subset P(\sigma)$, we have
\begin{equation}\label{uniform}
c \left( \frac{\ell_\sigma(A)}{\ell_\sigma(P(\sigma))} \right)< \frac{\ell_\tau(\T Q(A))}{\ell_\tau(\T Q (P(\sigma))} < \frac{1}{c}\left(\frac{\ell_\sigma(A)}{\ell_\sigma(P(\sigma))}\right).
\end{equation}

%%%%%%%%%%%%%%%%%%%%%%%%%%%%%%%%%%%%%%%%%%%%
\subsection{Proof of Theorem~\ref{limit-set}:}
%%%%%%%%%%%%%%%%%%%%%%%%%%%%%%%%%%%%%%%%%%%%
Cover $\pmf(S)$ by a finite collection of complete non-classical exchanges. Let $\tau$ be an exchange in this collection. By Theorem~\ref{proportion}, there is a subset $P(\tau, H) \subset P(\tau)$ that is disjoint from $\mathcal{R}_H$ and with proportion at least $K$. Consider the complement $P(\tau) \setminus P(\tau, H)$. By Theorem~\ref{uniform-distortion}, up to leaving out a set of $\ell_\tau$-measure zero, the complement $P(\tau) \setminus P(\tau, H)$ can be covered by $C$-uniformly distorted stages. Typically, there is an infinite number of these stages and we index these as $\tau_\alpha$. By applying Theorem~\ref{proportion} to each $\tau_\alpha$, we find subsets $P(\tau_\alpha, H) \subset P(\tau_\alpha)$ disjoint from $\mathcal{R}_H$ and with proportion at least $K$ in each $P(\tau_\alpha)$, where the proportion is being calculated in the measure $\ell_{\tau_\alpha}$. However, because $\tau_\alpha$ is $C$-uniformly distorted, estimate~\ref{uniform} implies that the proportion in the measure $\ell_\tau$ is at least $cK$. Then we split the complements $P(\tau_\alpha) \setminus P(\tau_\alpha, H)$ into $C$-uniformly distorted stages, and repeat the argument to find proportion $cK$ of these disjoint from $\mathcal{R}_H$, and so on. This completes the proof of Theorem~\ref{limit-set}.


\begin{thebibliography}{99}

\bibitem{Cal} Calegari, D. {\em Hyperbolic 3-manifolds, tameness, and Ahlfors' measure conjecture}, Lecture notes from a minicourse at the IUM, (2007).

\bibitem{Gad} Gadre, V. {\em Dynamics of non-classical interval exchanges}, http://arxiv.org/abs/0906.2563. 

\bibitem{Ker1} Kerckhoff, S. {\em Simplicial systems for interval exchange maps and measured foliations}, Ergodic Theory Dynam. Systems (1985), 5, 257-271.

\bibitem{Ker2} Kerckhoff, S. {\em The measure of the limit set of the handlebody group}, Topology, 29, (1990), no. 1, 27-40. 

\bibitem{Mas} Masur, H. {\em Measured foliations and handlebodies}, Ergodic Theory Dynam. Systems 6 (1986), no.1, 99-116. 

\bibitem{Pen} Penner, R. C. and Harer, J. L. {\em Combinatorics of train tracks}, Annals of Mathematics Studies, 125. Princeton University Press, Princeton, NJ, (1992).

\end{thebibliography}
\end{document}